\documentclass[11pt]{amsart}
\usepackage{amsfonts}
\usepackage{mathrsfs}
\usepackage{graphicx}
\usepackage{booktabs}
\usepackage{array}
\usepackage{amsmath}\allowdisplaybreaks[4]
\usepackage [latin1]{inputenc}
\usepackage[numbers,sort&compress]{natbib}

\newtheorem{thm}{Theorem}[section]
\newtheorem{cor}[thm]{Corollary}
\newtheorem{lem}[thm]{Lemma}
\newtheorem{prop}[thm]{Proposition}
\theoremstyle{definition}

\theoremstyle{remark}
\newtheorem{rem}[thm]{Remark}

\numberwithin{equation}{section}

\begin{document}

\title[Global solvability for the Boussinesq system]{{Global solvability for the Boussinesq system with fractional Laplacian}}
\author[Huiyang Zhang, Shuokai Yan and Qinghua Zhang]
{Huiyang Zhang, Shuokai Yan and Qinghua Zhang*} 

\thanks{{\it Address of the Authors:} School of Sciences,
  Nantong University,
6 Seyuan Road, Nantong City 226019, Jiangsu Province, P. R. China.}
\thanks{{\it * Corresponding author} E-mail: zhangqh@ntu.edu.cn}


\begin{abstract}
This paper focuses on the global solvability for the Boussinesq system with fractional Laplacian $(-\Delta)^{\alpha}$ in $\mathbb{R}^{n}$ for $n\geq3$. It proves the existence of a small positive number $\varepsilon=\varepsilon(n,\alpha)$ such that for each $0<T<\infty$, if $\frac{1}{2}<\alpha<\frac{2+n}{4}$ and $\|u_{0}\|_{\dot{H}^{s_{0}}}+T^{1/2}\|\theta_{0}\|_{\dot{H}^{s_{0}-\alpha}}\leq \varepsilon$, then the fractional Boussinesq system has a unique strong solution on the bounded interval $[0,T]$. If $\frac{1}{2}<\alpha<\frac{2+n}{6}$ and $\|u_{0}\|_{\dot{H}^{s_{0}}}+\|\theta_{0}\|_{\dot{H}^{s_{0}-2\alpha}}\leq \varepsilon$, then  the fractional Boussinesq system has a unique strong solution on the whole interval $[0,\infty)$.
\end{abstract}
\maketitle
\section{Introduction}
This paper aims at global solvability for  the Boussinesq system with fractional power of the negative Laplacian in $\mathbb{R}^{n}$ ($n\geq3$), i.e.
\begin{equation}\label{eqn:bsn}
\left\{\begin{array}{l}
\partial_{t}u+\mu(-\Delta)^{\alpha}u+u\cdot\nabla u+\nabla \pi=\kappa\theta e_{n},\;\;t>0,\;x\in\mathbb{R}^{n};\\
\partial_{t}\theta+\nu(-\Delta)^{\beta} \theta+u\cdot\nabla\theta=0,\;\;t>0,\;x\in\mathbb{R}^{n};\\
\nabla\cdot u=0,\;\;t>0,\;x\in\mathbb{R}^{n};\\
u(0,x)=u_{0}(x),\;\theta(0,x)=\theta_{0}(x),\;\;x\in\mathbb{R}^{n}.
\end{array} \right.
\end{equation}
Here $u(x,t)=(u_{1}(x,t),u_{2}(x,t),\cdots,u_{n}(x,t))$ denotes the velocity of the flow of a fluid at the place $x\in \mathbb{R}^{n}$  and the time $t>0$, while  $\theta(x,t)$ and $\pi(x,t)$ represent the temperature (or salinity) and inner pressure of the fluid. By $\theta_{0}(x)$ and $u_{0}(x)$ we denote the initial temperature (or salinity) and initial velocity respectively.

There are three physical constants related to the fluid: the proportion coefficient $\kappa$, which arises in the vortex buoyancy force $\kappa\theta e_{n}$, where $e_{n}=(0,\cdots,0,1)$ is the unit vertical vector, the viscous index $\mu$ and the thermal (or salinity ) confusion number $\nu$. For the sake of simplicity, we always assume that the three constants are all equal to $1$.

Here condition $\alpha=\beta$ is assumed, since in this case, Boussinesq system (\ref{eqn:bsn}) is variant under the scaling transformations
\begin{align*}
  u_{\lambda}(x,t)=\lambda^{2\alpha-1}u(\lambda x,\lambda^{2\alpha}t),\;\theta_{\lambda}(x,t)=\lambda^{4\alpha-1}\theta(\lambda x,\lambda^{2\alpha}t).
\end{align*}
Given two function spaces $A$ and $B$, we say $(A,B)$ is the critical space couple associated with the initial data of (\ref{eqn:bsn}), if the following equivalence relations hold
 \begin{align*}
   \|(u_{0})_{\lambda}\|_{A}\approx \|u_{0}\|_{A},\;\|(\theta_{0})_{\lambda}\|_{B}\approx\|\theta_{0}\|_{B}
 \end{align*}
There are two critical space couples for example, $(\dot{H}^{1-2\alpha+n/2}(\mathbb{R}^{n}),\dot{H}^{1-4\alpha+n/2}(\mathbb{R}^{n}))$  and $(\dot{B}_{p,r}^{1-2\alpha+n/p}(\mathbb{R}^{n}),\dot{B}_{p,r}^{1-4\alpha+n/p}(\mathbb{R}^{n}))$.

When $\theta=0$, system (\ref{eqn:bsn}) reduces to the Navier-Stokes equations with fractional Laplacian $(-\Delta)^{\alpha}$, or (GNS) in symbol. The well-known result is that (cf. \cite{li1969}) when $\alpha\geq\frac{n+2}{4}$, (GNS) has a unique global regular solution.
As for $\alpha<\frac{n+2}{4}$, the well-posedness for (GNS) reserves only for small data (small norm or small existing time), see \cite{wu2005,ch2017,du2018}, where global existence and uniqueness of the solutions with small norms in critical homogeneous Besov or Sobolev spaces were investigated.

Return to the Boussinesq systems. Global solvability of the 2D Boussinesq system was addressed in \cite{dp2011,mx2011,ah2007}. Global well-posedness of the 3D Boussinesq system with partial viscosity ($\nu=0$)  in critical Besov or Lorenz spaces was investigated in \cite{dp20081,dp20082} where $\alpha=1$, and in \cite{jy2016} where $\alpha\geq\frac{5}{4}$. For studies on the global solution to the 3D Boussinesq system with axisymmetric data, we refer to \cite{tr2010,mz2014} etc.

This paper focuses on the global solvability in homogeneous Sobolev spaces for the fractional Boussinesq system (\ref{eqn:bsn}) with full viscosity and $\beta=\alpha$ in $\mathbb{R}^{n}$ ($n\geq3$). By using the characterization of the space $\dot{H}^{s+\alpha}(\mathbb{R}^{n})$ and maximal $L^{2}-$regularity of the operator $(-\Delta)^{\alpha}$ on $\dot{H}^{s}(\mathbb{R}^{n})$ for $s<\frac{d}{2}$, traditional method of taking $L^{2}-$products in the equations to get an a prior estimates for the solutions is given up here.  We will prove that there exists  a small positive number $\varepsilon=\varepsilon(n,\alpha)$ such that for each $0<T<\infty$, if $\frac{1}{2}<\alpha<\frac{2+n}{4}$ and $\|u_{0}\|_{\dot{H}^{s_{0}}}+T^{1/2}\|\theta_{0}\|_{\dot{H}^{s_{0}-\alpha}}\leq \varepsilon$, then the initial value problem (\ref{eqn:bsn}) has a unique strong solution on the bounded interval $[0,T]$ (see Theorem \ref{thm:1.1}). If $\frac{1}{2}<\alpha<\frac{2+n}{6}$ and $\|u_{0}\|_{\dot{H}^{s_{0}}}+\|\theta_{0}\|_{\dot{H}^{s_{0}-2\alpha}}\leq \varepsilon$, then problem (\ref{eqn:bsn}) has a unique strong solution on the whole interval $[0,\infty)$  (see Theorem \ref{thm:1.2}).

\section{Preliminaries on function spaces and fractional Laplacian}

Let $\mathcal{S}$ be the collection of all rapidly decreasing smooth functions, and $\mathcal{S}'$ be its dual space, whose members are tempered generalised functions. Given a real number $s$, the homogeneous Sobolev space $\dot{H}^{s}:=\dot{H}^{s}(\mathbb{R}^{n})$ is defined by $\{f\in \mathcal{S}':\mathcal{F}^{-1}|\xi|^{s}\mathcal{F}\in L^{2}\}$, where $\mathcal{F}$ and $\mathcal{F}^{-1}$ denotes the $n-$dimensional Fourier transformation and its inverse respectively. For the sake of convenience, here we always use a single sign to denote both the scalar and vector function spaces, and omit the part $(\mathbb{R}^{n})$ in their notations.

Suppose $|s|<\frac{n}{2}$, and $p=\frac{2n}{n-2s}$, then we have the following embedding relations (cf. \cite[\S 1.3]{bcd2011})
 \begin{equation}\label{emb:1.1}
   \dot{H}^{s}\hookrightarrow L^{p}\;\mathrm{ if}\; 0\leq s<\frac{n}{2},\;\; L^{p}\hookrightarrow\dot{H}^{s}\;\mathrm{if}\; -\frac{n}{2}<s\leq 0.
\end{equation}
Moreover, if $s<\frac{n}{2}$, then $\dot{H}^{s}$ has a dense subset in $\mathcal{S}$, that is  $\mathcal{S}_{0}$,  whose members have the Fourier transformations vanishing near 0.

Define $\Lambda^{s}=\mathcal{F}^{-1}|\xi|^{s}\mathcal{F}$ and $(-\Delta)^{\alpha}=\Lambda^{2\alpha}$. Given $f\in \mathcal{S}'$, it is easy to see that $f$ lies in $\dot{H}^{s}(\mathbb{R}^{n})$ if and only if $\Lambda^{s}f\in L^{2}$. Furthermore, $\|(-\Delta)^{\alpha}f\|_{\dot{H}^{s}}=\|f\|_{\dot{H}^{s+2\alpha}}$, and for all $k\in \mathbb{N}$, the norms $\|\Lambda^{k}f\|_{2}$ and $\|\nabla^{k}f\|_{2}$ are equivalent. Define a commutator as follows
\begin{displaymath}
R_{s}(f,g)=\Lambda^{s}(fg)-(\Lambda^{s}f)g-f(\Lambda^{s}g),\;f,g\in \mathcal{S}_{0}.
\end{displaymath}
\begin{lem}[\cite{kpv1993}]
For any exponents $0<s,s_{1},s_{2}<1$, $1<p,q,r<+\infty$ verifying $s_{1}+s_{2}=1$ and $\frac{1}{p}=\frac{1}{q}+\frac{1}{r}$, there is corresponding a constant $C>0$ such that
\begin{equation}\label{ine:1.1}
 \|R_{s}(f,g)\|_{p}\leq C\|\Lambda^{s_{1}}f\|_{q}\|\Lambda^{s_{2}}g\|_{r}.
\end{equation}
\end{lem}
The following lemma is a corollary of \cite[Corollary 2.55]{bcd2011}. For the sake of completeness, we give a proof here without use of Besov spaces.
\begin{lem}
Assume that $|s|<\frac{n}{2}$, $\max\{s,0\}<s_{1},s_{2}<\frac{n}{2}$ such that $s_{1}+s_{2}=s+\frac{n}{2}$. Then for all $f\in\dot{H}^{s_{1}}(\mathbb{R}^{n})$,  $g\in\dot{H}^{s_{2}}(\mathbb{R}^{n})$, we have $fg\in \dot{H}^{s_{2}}(\mathbb{R}^{n})$ and
\begin{equation}\label{ine:product}
  \|fg\|_{\dot{H}^{s}(\mathbb{R}^{n})}\leq C\|f\|_{\dot{H}^{s_{1}}(\mathbb{R}^{n})}\|g\|_{\dot{H}^{s_{2}}(\mathbb{R}^{n})},
\end{equation}
where the constant $C>0$ is independent of $f$ and $g$.
\end{lem}
 \noindent{\bf Proof} : We will divide the proof into two cases.

 \noindent{\it Case 1.} $-\frac{n}{2}<s\leq0$.

 By the embedding inequality (\ref{ine:1.1}), it comes
 \begin{align*}
   \|fg\|_{\dot{H}^{s}(\mathbb{R}^{n})}\leq C\|fg\|_{\frac{2n}{n-2s}}\leq C\|f\|_{\frac{2n}{n-2s_{1}}}\|g\|_{\frac{2n}{n-2s_{2}}}\leq C\|f\|_{\dot{H}^{s_{1}}(\mathbb{R}^{n})}\|g\|_{\dot{H}^{s_{2}}(\mathbb{R}^{n})}.
 \end{align*}

 \noindent{\it Case 2.} $0<s<n/2$.

If $s\in \mathbb{N}$, then we have
\begin{align*}
\|fg\|_{\dot{H}^{s}}&\leq C\|\nabla^{s}(fg)\|_{2}\leq C\sum_{j=0}^{s}\|(\nabla^{j}f)(\nabla^{s-j}g)\|_{2}
\leq C\sum_{j=0}^{s}\|\nabla^{j}f\|_{\frac{2}{1-\alpha_{j}}}\|\nabla^{s-j}g\|_{\frac{2}{\alpha_{j}}}
\\&\leq C\sum_{j=0}^{s}\|\nabla^{j}f\|_{\dot{H}^{\frac{n}{2}\alpha_{j}}}\|\nabla^{s-j}g\|_{\dot{H}^{\frac{n}{2}(1-\alpha_{j})}}\leq C\|f\|_{\dot{H}^{s_{1}}}\|g\|_{\dot{H}^{s_{2}}},
\end{align*}
where $\alpha_{j}=\frac{2(s_{1}-j)}{n}$.

If $s\notin \mathbb{N}$, then let $k=[s]$ and $\delta=s-k$. By using (\ref{ine:1.1}), we can deduce that
\begin{align*}
\|fg\|_{\dot{H}^{s}}&\leq\|\nabla^{k}(fg)\|_{\dot{H}^{\delta}}\leq\sum_{j=0}^{k}\|(\nabla^{j}f)(\nabla^{k-j}g)\|_{\dot{H}^{\delta}}\\
&\leq\sum_{j=0}^{k}\big[\|R_{\delta}(\nabla^{j}f,\nabla^{k-j}g)\|_{2}+\|(\Lambda^{\delta}\nabla^{j}f)(\nabla^{k-j}g)\|_{2}+\|(\nabla^{j}f)(\Lambda^{\delta}\nabla^{k-j}g)\|_{2}\big]\\
&\leq C\sum_{j=0}^{k}\big[\|\Lambda^{\delta\theta_{j}}\nabla^{j}f\|_{\frac{2}{1-\theta_{j}}}\|\Lambda^{\delta(1-\theta_{j})}\nabla^{k-j}g\|_{\frac{2}{\theta_{j}}}\\
&\qquad+\|\Lambda^{\delta}\nabla^{j}f\|_{\frac{2}{1-\beta_{j}}}\|\nabla^{k-j}g\|_{\frac{2}{\beta_{j}}}+\|\nabla^{j}f\|_{\frac{2}{1-\alpha_{j}}}\|\Lambda^{\delta}\nabla^{k-j}g\|_{\frac{2}{\alpha_{j}}}\big]\\
&\leq C\sum_{j=0}^{k}\big[\|\Lambda^{\delta\theta_{j}}\nabla^{j}f\|_{\dot{H}^{\frac{n}{2}\theta_{j}}}\|\Lambda^{\delta(1-\theta_{j})}\nabla^{k-j}g\|_{\dot{H}^{\frac{n}{2}(1-\theta_{j})}}\\
&\qquad+C\|\Lambda^{\delta}\nabla^{j}f\|_{\dot{H}^{\frac{n}{2}\beta_{j}}}\|\nabla^{k-j}g\|_{\dot{H}^{\frac{n}{2}(1-\beta_{j})}}+C\|\nabla^{j}f\|_{\dot{H}^{\frac{n}{2}\alpha_{j}}}\|\Lambda^{\delta}\nabla^{k-j}g\|_{\dot{H}^{\frac{n}{2}(1-\alpha_{j})}}\big]\\
&\leq C\|f\|_{\dot{H}^{s_{1}}}\|g\|_{\dot{H}^{s_{2}}},
\end{align*}
where $\theta_{j}=\frac{2(s_{1}-j)}{n+2\delta}$ and $\beta_{j}=\frac{2(s_{1}-\delta-j)}{n}$.   \hfill $\Box$

 \begin{cor}
 Assume that $\frac{1}{2}<\alpha<\frac{1}{2}+\frac{n}{4}$ and let $s_{0}=1+\frac{n}{2}-2\alpha$. Then for a vector field $u$ and a scalar field $f$,  we can deduce that
\begin{equation}\label{ine:uf1}
  \|u\cdot\nabla f\|_{\dot{H}^{s_{0}-\alpha}}\leq C\|u\|_{\dot{H}^{s_{0}+\varepsilon}}\|\nabla f\|_{\dot{H}^{s_{0}+\alpha-1-\varepsilon}}\leq C\|u\|_{\dot{H}^{s_{0}+\varepsilon}}\|f\|_{\dot{H}^{s_{0}+\alpha-\varepsilon}}.
\end{equation}
where $0\leq\varepsilon<\min\{2\alpha-1,\alpha\}$ and the constant $C>0$ depends on $\varepsilon$.

Furthermore, if $\mathrm{div} u=0$, then for the same number $\varepsilon$, it holds
\begin{equation}\label{ine:uf2}
  \|u\cdot\nabla f\|_{\dot{H}^{s_{0}-2\alpha}}\leq C\|fu\|_{\dot{H}^{s_{0}-2\alpha+1}}\leq C\|u\|_{\dot{H}^{s_{0}+\varepsilon}}\|f\|_{\dot{H}^{s_{0}-\varepsilon}}.
\end{equation}
Assume further $\frac{1}{2}<\alpha<\frac{1}{3}+\frac{n}{6}$, then (\ref{ine:uf2}) turns to be
\begin{equation}\label{ine:uf3}
  \|u\cdot\nabla f\|_{\dot{H}^{s_{0}-3\alpha}}\leq C\|fu\|_{\dot{H}^{s_{0}-3\alpha+1}}\leq C\|u\|_{\dot{H}^{s_{0}+\varepsilon}}\|f\|_{\dot{H}^{s_{0}-\alpha-\varepsilon}}.
\end{equation}
\end{cor}

Let $s<\frac{n}{2}$, then the family of operators $\{e^{-t(-\Delta)^{\alpha}}=\mathcal{F}^{-1}e^{-t|\xi|^{2\alpha}}\mathcal{F},\;t\geq0\}$ defines a $C_{0}-$semigroup on $\dot{H}^{s}$, whose infinitesimal generator is $-(-\Delta)^{\alpha}$ with the domain $D((-\Delta)^{\alpha})=\dot{H}^{s}\cap\dot{H}^{s+2\alpha}$, where $s+2\alpha$ may be no less than $\frac{n}{2}$. Especially, if $s=0$, then $\dot{H}^{0}=L^{2}$ and $D((-\Delta)^{\alpha})=H^{2\alpha}$. As for the estimates of the semigroup $\{e^{-t(-\Delta)^{\alpha}}\}$, we have
\begin{align*}
\|(-\Delta)^{\gamma}e^{-t(-\Delta)^{\alpha}}f\|_{\dot{H}^{s}}=\||\xi|^{s+2\gamma}e^{-t|\xi|^{2\alpha}}\hat{f}\|_{2}\leq Ct^{-\gamma/\alpha}\|f\|_{\dot{H}^{s}},\;\;\forall\;t>0,
\end{align*}
where $\gamma\geq 0$ and the constant $C\geq1$ depends on $\gamma$. This shows that $\{e^{-t(-\Delta)^{\alpha}}\}$ is a uniformly bounded and analytic semigroup on $\dot{H}^{s}$.

For the sake of convenience, we denote $A_{\alpha}=(-\Delta)^{\alpha}$ in the following Lemma.
  \begin{lem}\label{lem:1.2}
$A_{\alpha}$ is a BIP type operator in $\dot{H}^{s}$, i.e. $A_{\alpha}^{iy}\in \mathcal{L}(\dot{H}^{s})$ for all $y\in\mathbb{R}$. Moreover,
\begin{equation}\label{bds:aiy}
  \|A_{\alpha}^{iy}\|_{\mathcal{L}(\dot{H}^{s})}\leq 1,\;\forall y\in \mathbb{R}.
\end{equation}
\end{lem}
\noindent{\bf Proof} : For each $\lambda>0$, the sum $A_{\alpha,\eta}=\eta I+A_{\alpha}$ is a positive operator, whose complex power $A_{\alpha,\eta}^{z}$ with $\mathrm{Re}z<0$ is defined below
\begin{equation*}
  A_{\alpha,\eta}^{z}=\frac{1}{2\pi i}\int_{\gamma}\lambda^{z}(\lambda I-A_{\alpha,\eta})^{-1}d\lambda,
\end{equation*}
where $\gamma=\gamma_{1}+\gamma_{2}+\gamma_{3}$, $\gamma_{1}$ and $\gamma_{3}$ are both half-lines from $\infty e^{\kappa i}$ to $\delta e^{\kappa i}$ and from $\delta e^{-\kappa i}$ to $\infty e^{-\kappa i}$ respectively, while $\gamma_{2}$ is a circular arc with the radius $\delta$ from $\delta e^{-\kappa\mathrm{i}}$ to $\delta e^{\kappa i}$, where $\pi/2<\kappa<\pi$ and $0<\delta<1$ such that the spectrum of $A_{\alpha,\eta}$ is located on the right side of $\gamma$. Furthermore, the function $\lambda^{z}$ takes the prime-value.

Given $f\in L^{2}$, we have
\begin{equation*}
  \mathcal{F}(A_{\alpha,\eta}^{z}f)=\frac{1}{2\pi i}\int_{\gamma}\lambda^{z}(\lambda +\eta+|\xi|^{2\alpha})^{-1}\mathcal{F}fd\lambda=(\eta+|\xi|^{2\alpha})^{z}\mathcal{F}f,
\end{equation*}
thus
\begin{equation*}
 A_{\alpha,\eta}^{z}f= \mathcal{F}^{-1}[(\eta+|\xi|^{2\alpha})^{z}\mathcal{F}f].
\end{equation*}
 As a straight consequence, it follows that
\begin{equation*}
  A_{\alpha,\eta}^{iy}=A_{\alpha,\eta}A_{\alpha,\eta}^{-1+iy}= \mathcal{F}^{-1}[(\eta+|\xi|^{2\alpha})^{it}\mathcal{F}f].
\end{equation*}
Since $|(\eta+|\xi|^{2\alpha})^{it}|\leq1$ for all $\eta>0$, $\xi\in\mathbb{R}^{n}$ and $t\in\mathbb{R}$, we see that $ A_{\alpha,\eta}^{iy}$ is an $L^{2}-$multiplier, and $\|A_{\alpha,\eta}^{iy}\|_{\mathcal{L}(L^{2})}\leq 1$ for all $y\in \mathbb{R}$, which in turn yields
\begin{equation}\label{est:aiy}
  \|A_{\alpha,\eta}^{iy}\|_{\mathcal{L}(\dot{H}^{s})}=\|A_{\alpha,\eta}^{iy}\|_{\mathcal{L}(L^{2})}\leq 1,\;\forall y\in \mathbb{R}.
\end{equation}

Note that $A_{\alpha}$ is an injective operator with $D(A_{\alpha})$ and $R(A_{\alpha})$ both dense in $\dot{H}^{s}$, from \cite[Theorem 7.64 \& Proposition 7.47]{ms2001}, we can deduce that $A_{\alpha}^{iy}\in\mathcal{L}(\dot{H}^{s})$, and
$A_{\alpha}^{iy}=s-\lim_{\eta\rightarrow0^{+}}A_{\alpha,\eta}^{iy}$, which combined with (\ref{est:aiy}), leads to the desired estimates (\ref{bds:aiy}).  \hfill $\Box$

Now by invoking \cite{gs1991}, we conclude that $(-\Delta)^{\alpha}$ has the $L^{2}-$maximal regularity on $[0,T]$ for all $0<T\leq\infty$, i.e.
 \begin{prop}\label{prop:2.1}
For all $f\in L^{2}(0,T;\dot{H}^{s})$, there is a unique solution $u$ to the following Cauchy problem
\begin{equation}\label{eqn:lin0}
         \partial_{t}w+(-\Delta)^{\alpha}w=f(t)\;\textrm{for a.e.}\;t>0, \;\;
        w(0)=0
\end{equation}
such that $w\in L(0,T;\dot{H}^{s+2\alpha})$, $w'\in L^{2}(0,T;\dot{H}^{s})$ and
\begin{equation}\label{ine:mr2}
 \|\partial_{t}w\|_{L^{2}(0,T;\dot{H}^{s})}+\|w\|_{L^{2}(0,T;\dot{H}^{s+2\alpha})}\leq C\|f\|_{L^{2}(0,T;\dot{H}^{s})},
\end{equation}
where $\partial_{t}w$ is the derivative of $w$ w.r.t. t in the sense of distribution, $\|w\|_{L^{2}(0,T;\dot{H}^{s+2\alpha})}=\|(-\Delta)^{\alpha}w\|_{L^{2}(0,T;\dot{H}^{s})}$, and the constant $C>0$ is independent of $T$ and $f$.
\end{prop}
 \begin{rem}
 By means of Fourier transformation, one can easily check that the unique solution to (\ref{eqn:lin0}) can be represented by
 \begin{equation*}
   w(t)=\int_{0}^{t}e^{-(t-\tau)(-\Delta)^{\alpha}}f(\tau)d\tau.
 \end{equation*}
\end{rem}

Let $s+\alpha<\frac{n}{2}$, then for $w\in L^{2}(0,T;\dot{H}^{s+2\alpha})$ with $\partial_{t}w\in L^{2}(0,T;\dot{H}^{s})$, we have
\begin{equation*}
 \frac{d}{dt}\|w(t)\|_{\dot{H}^{s+\alpha}}^{2}=2(\Lambda^{s}\partial_{t}w,\Lambda^{s+2\alpha}w(t))\;\textrm{for a.e.}\;t\in(0,T),
\end{equation*}
from which we can deduce that if $w(0)\in \dot{H}^{s+\alpha}$, then $w$ is bounded and uniformly continuous on $[0,T)$ in $\dot{H}^{s+\alpha}$, or $w\in BUC([0,T);\dot{H}^{s+\alpha})$ in symbol, and
\begin{equation}\label{est:wt}
  \sup_{t\in[0,T)}\|w(t)\|_{\dot{H}^{s+\alpha}}^{2}\leq \big(\|w(0)\|_{\dot{H}^{s+\alpha}}^{2}+2\|w\|_{L^{2}(0,T;\dot{H}^{s+2\alpha})}\|\partial_{t}w\|_{L^{2}(0,T;\dot{H}^{s})}\big).
\end{equation}

Direct calculation also shows that $a\in\dot{H}^{s+\alpha}$  if and only if $\|(-\Delta)^{\alpha}e^{-t(-\Delta)^{\alpha}}a\|_{\dot{H}^{s}}\in L^{2}(0,\infty)$, and
\begin{equation}\label{ch:int}
 \|a\|_{\dot{H}^{s+\alpha}}^{2}=\frac{1}{2}\int_{0}^{\infty}\|(-\Delta)^{\alpha}e^{-t(-\Delta)^{\alpha}}a\|_{\dot{H}^{s}}^{2}dt.
\end{equation}
Therefore by the analyticity of $e^{-t(-\Delta)^{\alpha}}$, we can deduce that
\begin{prop}\label{prop:2.2}
For each $a\in\dot{H}^{s+\alpha}$, the function $a_{L}:=e^{-t(-\Delta)^{\alpha}}a$  is the strong solution to the problem
 \begin{equation*}
         \partial_{t}v+(-\Delta)^{\alpha}v=0,\;\forall\;t>0, \;\;
       v(0)=a.
\end{equation*}
Moreover, we have $a_{L}\in BUC([0,\infty);H^{s+\alpha})\cap L^{2}(0,\infty;\dot{H}^{s+2\alpha})$ such that $\partial_{t}a_{L}\in L^{2}(0,\infty;\dot{H}^{s})$, and
\begin{equation}\label{est:al}
\sup_{t\in[0,\infty)}\|a_{L}(t)\|_{ \dot{H}^{s+\alpha}}+\|a_{L}\|_{L^{2}(0,\infty;\dot{H}^{s+2\alpha})}+\|\partial_{t}a_{L}\|_{L^{2}(0,\infty;\dot{H}^{s})}\leq 3\|a\|_{\dot{H}^{s+\alpha}}.
\end{equation}
\end{prop}

\section{Existence and uniqueness of global strong solution}

Given $s<n/2$, denote $\dot{H}_{\sigma}^{s}$ by the collection of all solenoidal vector fields whose components lie in $\dot{H}^{s}$. Introduce two operators $Q$ and $P$ on the space $\dot{H}^{s}$ of vector type as follows
\begin{align*}
Qu=\nabla(-\Delta)^{-1}(\nabla\cdot u)=\mathcal{F}^{-1}\big(\frac{i\xi}{|\xi|^{2}}\xi\cdot \mathcal{F}u\big)\;\mathrm{and}\;Pu=u-Qu.
\end{align*}
Note that $|\xi_{j}\xi_{k}|/|\xi|^{2}\leq1$ for all $\xi\in\mathbb{R}^{n}\setminus\{0\}$ and all $k,j=i,2,\cdot,n$, we see that $P,Q$ are both $L^{2}-$ and consequently $\dot{H}^{s}-$multipliers, and $\|Qu\|_{\dot{H}^{s}}\leq \|u\|_{\dot{H}^{s}}$, $\|Pu\|_{\dot{H}^{s}}\leq \|u\|_{\dot{H}^{s}}$. We can also check that, $Pu\in \dot{H}_{\sigma}^{s}$ and $(Pu,Qu)_{\dot{H}^{s}}=0$ for all $u\in\dot{H}^{s}$, and  $P^{2}=P$, $Q^{2}=Q$. Based on these properties,  we can derive the equality $\dot{H}_{\sigma}^{s}=P\dot{H}^{s}$ and the Helmholtz decomposition: $\dot{H}_{\sigma}^{s}=\dot{H}_{\sigma}^{s}\oplus Q\dot{H}^{s}$. Here $P:\dot{H}^{s}\rightarrow\dot{H}_{\sigma}^{s}$ is call the  Helmholtz projection. It is commutative with $(-\Delta)^{\alpha}$ on $\mathbb{R}^{n}$.

Let $s_{i}<n/2$, $i=1,2$, $(u_{0},\theta_{0})\in \dot{H}_{\sigma}^{s_{1}}\times\dot{H}^{s_{2}}$ and $0<T\leq\infty$. We say the function triple $(u,\nabla\pi,\theta)$ is a strong solution to the initial value problem of Boussinesq system (\ref{eqn:bsn}) on $[0,T)$, if  $u\in C([0,T);\dot{H}_{\sigma}^{s_{1}})\cap L_{\mathrm{loc}}^{1}(0,T;\dot{H}^{s_{1}+\alpha})$, $\nabla\pi\in L_{\mathrm{loc}}^{1}(0,T;\dot{H}^{s_{1}-\alpha})$ and $\theta\in C([0,T);\dot{H}^{s_{2}})\cap L_{\mathrm{loc}}^{1}(0,T;\dot{H}^{s_{2}+\alpha})$ such that  $\partial_{t}u\in L_{\mathrm{loc}}^{1}(0,T;\dot{H}^{s_{1}-\alpha})$, $\partial_{t}\theta\in L_{\mathrm{loc}}^{1}(0,T;\dot{H}^{s_{2}-\alpha})$, and
\begin{equation}\label{eqn:bsn3}
\left\{\begin{array}{l}
\partial_{t}u+\mu(-\Delta)^{\alpha}u+(u\cdot\nabla u)+\nabla\pi=\theta e_{n}\;\;\mathrm{in}\;\;\dot{H}^{s_{1}-\alpha},\\
\partial_{t}\theta+\nu(-\Delta)^{\alpha} \theta+u\cdot\nabla\theta=0\;\;\mathrm{in}\;\;\dot{H}^{s_{2}-\alpha}
\end{array} \right.
\end{equation}
for a.e. $t\in[0,T)$, and
$u(0)=u_{0}$, $\theta(0)=\theta_{0}$.

 \begin{thm}\label{thm:1.1}
Suppose that $\frac{1}{2}<\alpha<\frac{2+n}{4}$. Then there is a small number $\varepsilon=\varepsilon(n,\alpha)>0$ such that for every $0<T<\infty$, if $(u_{0},\theta_{0})\in \dot{H}_{\sigma}^{s_{0}}\times\dot{H}^{s_{0}-\alpha}$ verifying
 \begin{equation}\label{cd:ini1}
   \|u_{0}\|_{\dot{H}^{s_{0}}}+T^{1/2}\|\theta_{0}\|_{\dot{H}^{s_{0}-\alpha}}\leq\varepsilon,
 \end{equation}
 then on the interval $[0,T]$, the initial value problem (\ref{eqn:bsn}) has a unique strong solution $(u,\nabla\pi,\theta)$, up to an additional constant to $\pi$. This solution lies in the space $X_{T}\times L^{2}(0,T;\dot{H}^{s_{0}-\alpha})\times Y_{T}$ and verifies the following estimates
 \begin{equation*}
   \|u\|_{X_{T}}+\|\nabla\pi\|_{ L^{2}(0,T;\dot{H}^{s_{0}-\alpha})}+T^{1/2}\|\theta\|_{Y_{T}}\leq C(\|u_{0}\|_{\dot{H}^{s_{0}}}+T^{1/2}\|\theta_{0}\|_{\dot{H}^{s_{0}-\alpha}})
 \end{equation*}
for some constant constant $C>0$ independent of $u_{0}$, $\theta_{0}$ and $T$.

 Furthermore, uniqueness of the strong solution to (\ref{eqn:bsn}) in $X_{T}\times L^{2}(0,T;\dot{H}^{s_{0}-\alpha})\times Y_{T}$ holds true without the initial assumption (\ref{cd:ini1}).
\end{thm}
\noindent{\bf Proof} :
Given $0<T<\infty$, we introduce two function spaces
  \begin{align*}
&X_{T}=\big\{u\in C([0,T];\dot{H}_{\sigma}^{s_{0}})\cap L^{2}(0,T;\dot{H}^{s_{0}+\alpha}):\partial_{t}u\in L^{2}(0,T;\dot{H}^{s_{0}-\alpha})\big\},\\
&Y_{T}=\big\{\theta\in  C([0,T];\dot{H}^{s_{0}-\alpha})\cap L^{2}(0,T;\dot{H}^{s_{0}}):\partial_{t}\theta\in L^{2}(0,T;\dot{H}^{s_{0}-2\alpha})\big\}
 \end{align*}
 with the norms defined by
 \begin{align*}
   &\|u\|_{X_{T}}=\max_{0\leq t\leq T}\|u(t)\|_{\dot{H}^{s_{0}}}+\|u\|_{L^{2}(0,T;\dot{H}^{s_{0}+\alpha})}+\|\partial_{t}u\|_{L^{2}(0,T;\dot{H}^{s_{0}-\alpha})},\\
   &\|\theta\|_{Y_{T}}=\max_{0\leq t\leq T}\|\theta(t)\|_{\dot{H}^{s_{0}-\alpha}}+\|\theta\|_{L^{2}(0,T;\dot{H}^{s_{0}})}+\|\partial_{t}\theta\|_{L^{2}(0,T;\dot{H}^{s_{0}-2\alpha})}.
 \end{align*}

This proof will be divided into two parts.

\noindent {\it (1) Existence part}.

We first search for the strong solution of the following abstract evolution equations
\begin{equation}\label{eqn:bsn1}
\left\{\begin{array}{l}
\partial_{t}u+\mu(-\Delta)^{\alpha}u+P(u\cdot\nabla u)= P(\theta e_{n})\;\;\mathrm{in}\;\;\dot{H}_{\sigma}^{s_{0}-\alpha},\;\mathrm{a.e.}\;t>0,\\
\partial_{t}\theta+\nu(-\Delta)^{\alpha} \theta+u\cdot\nabla\theta=0\;\;\mathrm{in}\;\;\dot{H}^{s_{0}-2\alpha},\;\mathrm{a.e.}\;t>0,\\
u(0)=u_{0},\;\theta(0)=\theta_{0}.
\end{array} \right.
\end{equation}
From Proposition \ref{prop:2.1}, \ref{prop:2.2}, we see that every solution of (\ref{eqn:bsn1}) also solves the  integral equations
 \begin{align}\label{eqn:bsn2}
u=u_{L}+L(\theta)+\Phi(u,u),\;\;\theta=\theta_{L}+\Psi(u,\theta),
 \end{align}
 where $u_{L}=e^{-t(-\Delta)^{\alpha}}u_{0}$, $\theta_{L}=e^{-t(-\Delta)^{\alpha}}\theta_{0}$, and
 \begin{align*}
&L(\theta)=\int_{0}^{t}e^{-(t-\tau)(-\Delta)^{\alpha}} P(\theta e_{n})d\tau,\\
&\Phi(u,v)=-\int_{0}^{t}e^{-(t-\tau)(-\Delta)^{\alpha}}P(u\cdot\nabla u)d\tau,\\
&\Psi(u,\theta)=-\int_{0}^{t}e^{-(t-\tau)(-\Delta)^{\alpha}} u\cdot\nabla \theta d\tau.
\end{align*}

By taking (\ref{est:al}) into account, we have $u_{L}\in X_{T}$, $\theta_{L}\in Y_{T}$ and
\begin{equation}\label{est:uthe1}
 \|u_{L}\|_{X_{T}}\leq 3\|u_{0}\|_{\dot{H}^{s_{0}}},\;\|\theta_{L}\|_{Y_{T}} \leq 3\|\theta_{0}\|_{\dot{H}^{s_{0}-\alpha}}.
\end{equation}
Besides, for all $u,v\in X_{T}$ and $\theta\in Y_{T}$, by employing the inequalities (\ref{ine:uf1}), (\ref{ine:uf2}) together with the estimate (\ref{ine:mr2}), we can deduce that $L(\theta)\in X_{T}$, $\Phi(u,v)\in X_{T}$ and $\Psi(u,\theta)\in Y_{T}$, and
\begin{align}\label{est:l}
\|L(\theta)\|_{X_{T}}&\leq C\|P(\theta e_{n})\|_{L^{2}(0,T;\dot{H}^{s_{0}-\alpha})}\leq k_{1}T^{1/2}\max_{0\leq t\leq T}\|\theta(t)\|_{\dot{H}^{s_{0}-\alpha}}\\
\nonumber&\leq k_{1}T^{1/2}\|\theta\|_{Y_{T}},
\end{align}
\begin{align}\label{est:phi}
\|\Phi(u,v)\|_{X_{T}}&\leq C\|P(u\cdot\nabla v)\|_{L^{2}(0,T;\dot{H}^{s_{0}-\alpha})}\\
\nonumber&\leq k_{2}\max_{0\leq t\leq T}\|u(t)\|_{\dot{H}^{s_{0}}}
\|v\|_{L^{2}(0,T;\dot{H}^{s_{0}+\alpha})}\leq k_{2}\|u\|_{X_{T}}\|v\|_{X_{T}},
\end{align}
\begin{align}\label{est:psi}
\|\Psi(u,\theta)\|_{Y_{T}}&\leq C\|u\cdot\nabla \theta\|_{L^{2}(0,T;\dot{H}^{s_{0}-2\alpha})}\\
\nonumber&\leq k_{3}\max_{0\leq t\leq T}\|u(t)\|_{\dot{H}^{s_{0}}}\|\theta\|_{L^{2}(0,T;\dot{H}^{s_{0}})}\leq k_{3}\|u\|_{X_{T}}\|\theta\|_{Y_{T}},
\end{align}
where the constants $k_{i}>0$, $i=1,2,3$ are all independent of $T$, $u$, $v$ and $\theta$.

Letting $\|\theta\|_{Y_{T}}'=2k_{1}T^{1/2}\|\theta\|_{Y_{T}}$, estimates (\ref{est:l}) and (\ref{est:psi}) turn to be
\begin{equation*}
\|L(\theta)\|_{X_{T}}\leq\frac{1}{2}\|\theta\|_{Y_{T}}'\;\textrm{and}\;\|\Psi(u,\theta)\|_{Y_{T}'}\leq k_{3}\|u\|_{X_{T}}\|\theta\|_{Y_{T}}'.
\end{equation*}
Thus if we assume
\begin{equation}\label{h:k0}
  K_{0}:=\|u_{L}\|_{X_{T}}+\|\theta_{L}\|_{Y_{T}}'\leq\frac{1}{32(k_{2}+k_{3})}
\end{equation}
and take
\begin{equation*}
  \lambda_{1}=\frac{1-\sqrt{1-16K_{0}(k_{2}+k_{3})}}{4(k_{2}+k_{3})},
\end{equation*}
the first root of the quadratic equation
$(k_{2}+k_{3})\lambda^{2}-\lambda/2+K_{0}=0$,
then by invoking \cite[Lemma 5.1]{bs2012}, we can deduce that the integral equations (\ref{eqn:bsn2})
have a unique solution $(u,\theta)\in X_{T}\times Y_{T}$ obeying the following estimate
\begin{equation}\label{est:uth}
  \|u\|_{X_{T}}+2k_{1}T^{1/2}\|\theta\|_{Y_{T}}\leq \lambda_{1}.
\end{equation}
Note that with neglect of the difference between scalar and vector-valued functions, $L(\theta)$, $\Phi(u,v)$ and $\Psi(u,\theta)$ are all solutions of (\ref{eqn:lin0}) with $f=P(\theta e_{n})$, $-P(u\cdot\nabla u)$ and $u\cdot\nabla \theta$ respectively.
Hence the function couple $(u,\theta)$ are exactly the strong solution of the Cauchy problem (\ref{eqn:bsn1}).

Now taking (\ref{est:uthe1}) into account, we conclude that under the assumption
\begin{equation*}
  \|u_{0}\|_{\dot{H}^{s_{0}}}+2k_{1}T^{1/2}\|\theta_{0}\|_{\dot{H}^{s_{0}-\alpha}}\leq\frac{1}{96(k_{2}+k_{3})},
\end{equation*}
hypothesis (\ref{h:k0}) is fulfilled, consequently (\ref{eqn:bsn1}) admits a strong solution $(u,\theta)$ in $X_{T}\times Y_{T}$. Additionally, using (\ref{est:uthe1}) again, estimate (\ref{est:uth}) becomes
 \begin{align}\label{est:sol}
 \|u\|_{X_{T}}+2k_{1}T^{1/2}\|\theta\|_{Y_{T}}\leq\lambda_{1}\leq4K_{0}\leq12(\|u_{0}\|_{\dot{H}^{s_{0}}}+2k_{1}T^{1/2}\|\theta_{0}\|_{\dot{H}^{s_{0}-\alpha}}).
 \end{align}

 Let $\nabla\pi=Q(\theta e_{n}-u\cdot\nabla u)$, then we have $\nabla\pi\in L^{2}(0,T;\dot{H}^{s_{0}-\alpha})$, and
 \begin{align*}
   \|\nabla\pi\|_{L^{2}(0,T;\dot{H}^{s_{0}-\alpha})}&\leq \|\theta e_{n}\|_{L^{2}(0,T;\dot{H}^{s_{0}-\alpha})}+\|u\cdot\nabla u\|_{L^{2}(0,T;\dot{H}^{s_{0}-\alpha})}\\
   &\leq \frac{1}{2}\|\theta\|_{Y_{T}}'+k_{2}\|u\|_{X_{T}}^{2}
   \leq\frac{\lambda_{1}}{2}+k_{2}\lambda_{1}^{2}\\
   &\leq\lambda_{1}\leq12(\|u_{0}\|_{\dot{H}^{s_{0}}}+2k_{1}T^{1/2}\|\theta_{0}\|_{\dot{H}^{s_{0}-\alpha}}).
 \end{align*}
 It is easy to check that the function triple $(u,\nabla\pi,\theta)$ verifies both the equations of (\ref{eqn:bsn3}) in $\dot{H}^{s_{0}-\alpha}$ and $\dot{H}^{s_{0}-2\alpha}$ for a.e. $t\in[0,T]$. Thus existence of the strong solution of (\ref{eqn:bsn}) has been reached.

\noindent {\it (2) Uniqueness part}.

Assume that  $(u_{1},\nabla\pi_{1},\theta_{1})$ and $(u_{2},\nabla\pi_{2},\theta_{2})$ are two solutions of the system (\ref{eqn:bsn}). Consider the differences $\delta u=u_{1}-u_{2}$, $\nabla\delta \pi=\nabla\pi_{1}-\nabla\pi_{2}$ and $\delta \theta=\theta_{1}-\theta_{2}$. Then the function triple $(\delta u,\nabla\delta \pi,\delta \theta)$ is the strong solution of the homogeneous system
\begin{align*}
\left\{\begin{array}{l}
\partial_{t}w+(-\Delta)^{\alpha}w+\nabla \varpi=\vartheta e_{n}-u_{1}\cdot\nabla w-w\cdot\nabla u_{2}\;\;\mathrm{in}\;\;\dot{H}^{s_{0}-\alpha},\;\mathrm{a.e.}\;t>0,\\
\partial_{t}\vartheta+\nu(-\Delta)^{\beta} \vartheta=-u_{1}\cdot\nabla\vartheta-w\cdot\nabla\theta_{2}\;\;\mathrm{in}\;\;\dot{H}^{s_{0}-2\alpha},\;\mathrm{a.e.}\;t>0,\\
w(0)=0,\;\vartheta(0,x)=0.
\end{array} \right.
\end{align*}
Consequently,  the function couple $(\delta u,\delta \theta)$ solves the following system
\begin{align}\label{eqn:dut}
w=L(\vartheta)+\Phi(u_{1},w)+\Phi(w,u_{2}),\;\;\vartheta=\Psi(u_{1},\vartheta)+\Psi(w,\vartheta_{2}).
 \end{align}

With the aid of (\ref{est:l}), (\ref{est:phi}) and (\ref{est:psi}), together with (\ref{ine:uf1}), we can derive from (\ref{eqn:dut}) that
 \begin{align}\label{est:du}
\nonumber\|\delta u\|_{X_{t}}&\leq k_{1}T^{1/2}\|\delta\theta\|_{Y_{t}}+k_{2}\|\delta u\|_{X_{t}}\|u_{2}\|_{L^{2}(0,t;\dot{H}^{s_{0}+\alpha})}\\
&\quad+C\big\|\|u_{1}\|_{\dot{H}^{s_{0}+\varepsilon}}\|\delta u\|_{\dot{H}^{s_{0}+\alpha-\varepsilon}}\big\|_{L^{2}(0,t)}
\end{align}
and
\begin{align}\label{est:dth}
\|\delta \theta\|_{Y_{t}}\leq k_{3}\|\delta u\|_{X_{t}}\|\theta_{2}\|_{L^{2}(0,t;\dot{H}^{s_{0}})}+C\big\|\|u_{1}\|_{\dot{H}^{s_{0}+\varepsilon}}\|\delta\theta\|_{\dot{H}^{s_{0}-\varepsilon}}\big\|_{L^{2}(0,t)}.
\end{align}
By means of interpolations, we have
 \begin{align*}
&\quad\big\|\|u_{1}\|_{\dot{H}^{s_{0}+\varepsilon}}\|\delta u\|_{\dot{H}^{s_{0}+\alpha-\varepsilon}}\big\|_{L^{2}(0,t)}\\
&\leq\big\|\|u_{1}\|_{\dot{H}^{s_{0}}}^{1-\varepsilon/\alpha}\|u_{1}\|_{\dot{H}^{s_{0}+\alpha}}^{\varepsilon/\alpha}\|\delta u\|_{\dot{H}^{s_{0}}}^{\varepsilon/\alpha}\|\delta u\|_{\dot{H}^{s_{0}+\alpha}}^{1-\varepsilon/\alpha}\big\|_{L^{2}(0,t)}\\
&\leq \max_{0\leq \tau\leq t}\|u_{1}(\tau)\|_{\dot{H}^{s_{0}}}^{1-\varepsilon/\alpha}\|\delta u(\tau)\|_{\dot{H}^{s_{0}}}^{\varepsilon/\alpha}\cdot\|u_{1}\|_{L^{2}(0,t;\dot{H}^{s_{0}+\alpha})}^{\varepsilon/\alpha}\|\delta u\|_{L^{2}(0,t;\dot{H}^{s_{0}+\alpha})}^{1-\varepsilon/\alpha}\\
&\leq\|u_{1}\|_{X_{T}}^{1-\varepsilon/\alpha}\|u_{1}\|_{L^{2}(0,t;\dot{H}^{s_{0}+\alpha})}^{\varepsilon/\alpha}\|\delta u\|_{X_{t}},
\end{align*}
and
 \begin{align*}
&\quad\big\|\|u_{1}\|_{\dot{H}^{s_{0}+\varepsilon}}\|\delta\theta\|_{\dot{H}^{s_{0}-\varepsilon}}\big\|_{L^{2}(0,t)}\\
&\leq\big\|\|u_{1}\|_{\dot{H}^{s_{0}}}^{1-\varepsilon/\alpha}\|u_{1}\|_{\dot{H}^{s_{0}+\alpha}}^{\varepsilon/\alpha}\|\delta \theta\|_{\dot{H}^{s_{0}-\alpha}}^{\varepsilon/\alpha}\|\delta\theta\|_{\dot{H}^{s_{0}}}^{1-\varepsilon/\alpha}\big\|_{L^{2}(0,t)}\\
&\leq \max_{0\leq \tau\leq t}\|u_{1}(\tau)\|_{\dot{H}^{s_{0}}}^{1-\varepsilon/\alpha}\|\delta \theta(\tau)\|_{\dot{H}^{s_{0}-\alpha}}^{\varepsilon/\alpha}\cdot\|u_{1}\|_{L^{2}(0,t;\dot{H}^{s_{0}+\alpha})}^{\varepsilon/\alpha}\|\delta \theta\|_{L^{2}(0,t;\dot{H}^{s_{0}})}^{1-\varepsilon/\alpha}\\
&\leq\|u_{1}\|_{X_{T}}^{1-\varepsilon/\alpha}\|u_{1}\|_{L^{2}(0,t;\dot{H}^{s_{0}+\alpha})}^{\varepsilon/\alpha}\|\delta \theta\|_{Y_{t}}.
\end{align*}
Since $\|u_{i}\|_{L^{2}(0,t;\dot{H}^{s_{0}+\alpha})}=o(1)$ and $\|\theta_{i}\|_{L^{2}(0,t;\dot{H}^{s_{0}})}=o(1)$ as $t\rightarrow0$, $i=1,2$, we can find a small number $t_{0}>0$ such that
\begin{align*}
k_{2}\|u_{2}\|_{L^{2}(0,t;\dot{H}^{s_{0}+\alpha})}+k_{3}\|\theta_{2}\|_{L^{2}(0,t;\dot{H}^{s_{0}})}+C\|u_{1}\|_{X_{T}}^{1-\varepsilon/\alpha}\|u_{1}\|_{L^{2}(0,t;\dot{H}^{s_{0}+\alpha})}^{\varepsilon/\alpha}\leq\frac{1}{2}
\end{align*}
 for all $0\leq t\leq t_{0}$. Putting them into (\ref{est:du}) and (\ref{est:dth}), we can deduce that
 \begin{align*}
\|\delta u\|_{X_{t_{0}}}=\|\delta\theta\|_{Y_{t_{0}}}=0.
\end{align*}

Using $t_{0}$ as the new initial time, and performing the same arguments, we can find another time $t_{1}>t_{0}$ such that
\begin{align*}
\|\delta u\|_{X_{t_{1}}}=\|\delta u(\cdot+t_{0})\|_{X_{t_{1}-t_{0}}}=0\;\;\mathrm{and}\;\;\|\delta \theta\|_{X_{t_{1}}}=\|\delta \theta(\cdot+t_{0})\|_{X_{t_{1}-t_{0}}}=0.
\end{align*}
Thus based on the connection of $[0,T]$ and continuity of the norms $\|u_{i}\|_{L^{2}(0,t;\dot{H}^{s_{0}+\alpha})}$ and $\|\theta_{i}\|_{L^{2}(0,t;\dot{H}^{s_{0}})}$ w.r.t. t, and by means of iteration, we can eventually prove that $\delta u=0$ in $X_{T}$ and $\delta\theta=0$ in $Y_{T}$. Additionally, since
\begin{align*}
\nabla\delta\pi=Q(\delta\theta e_{n}-u_{1}\cdot\nabla \delta u-\delta u\cdot\nabla u_{2})
\end{align*}
 we obtain $\nabla\delta\pi=0$ in $L^{2}(0,T;\dot{H}^{s_{0}-\alpha})$. Thus uniqueness of the strong solution has been reached.
\hfill $\Box$

 \begin{thm}\label{thm:1.2}
Suppose that $\frac{1}{2}<\alpha<\frac{1}{3}+\frac{n}{6}$. Then there is a small number $\varepsilon=\varepsilon(n,\alpha)>0$ such that if the initial data $(u_{0},\theta_{0})\in \dot{H}_{\sigma}^{s_{0}}\times\dot{H}^{s_{0}-2\alpha}$, and
 \begin{equation}\label{cd:ini2}
   \|u_{0}\|_{\dot{H}^{s_{0}}}+\|\theta_{0}\|_{\dot{H}^{s_{0}-2\alpha}}\leq\varepsilon,
 \end{equation}
the initial value problem (\ref{eqn:bsn}) has a unique global strong solution $(u,\nabla\pi,\theta)$ in the class
$X\times L^{2}(0,\infty;\dot{H}^{s_{0}-\alpha})\times Y$, and
\begin{align*}
 \|u\|_{X}+\|\nabla\pi\|_{ L^{2}(0,\infty;\dot{H}^{s_{0}-\alpha})}+\|\theta\|_{Y}\leq C(\|u_{0}\|_{\dot{H}^{s_{0}}}+\|\theta_{0}\|_{\dot{H}^{s_{0}-2\alpha}}),
 \end{align*}
 where the constant $C>0$ is independent of $u_{0}$ and $\theta_{0}$.

 Furthermore, uniqueness of the strong solution to (\ref{eqn:bsn}) in $X\times L^{2}(0,\infty;\dot{H}^{s_{0}-\alpha})\times Y$ holds true without the initial assumption (\ref{cd:ini2}).
\end{thm}
\noindent{\bf Proof} :
Introduce other two function spaces
  \begin{align*}
&X=\big\{u\in BUC([0,\infty);\dot{H}_{\sigma}^{s_{0}})\cap L^{2}(0,\infty;\dot{H}^{s_{0}+\alpha}):\partial_{t}u\in L^{2}(0,\infty;\dot{H}^{s_{0}-\alpha})\big\},\\
&Y=\big\{\theta\in  BUC([0,\infty);\dot{H}^{s_{0}-2\alpha})\cap L^{2}(0,\infty;\dot{H}^{s_{0}-\alpha}):\partial_{t}\theta\in L^{2}(0,\infty;\dot{H}^{s_{0}-3\alpha})\big\},
 \end{align*}
whose norms are defined by
 \begin{align*}
   &\|u\|_{X}=\sup_{t\geq0}\|u(t)\|_{\dot{H}^{s_{0}}}+\|u\|_{L^{2}(0,\infty;\dot{H}^{s_{0}+\alpha})}+\|\partial_{t}u\|_{L^{2}(0,\infty;\dot{H}^{s_{0}-\alpha})},\\
   &\|\theta\|_{Y}=\sup_{t\geq0}\|\theta\|_{\dot{H}^{s_{0}-2\alpha}}+\|\theta\|_{L^{2}(0,\infty;\dot{H}^{s_{0}-\alpha})}+\|\partial_{t}\theta\|_{L^{2}(0,\infty;\dot{H}^{s_{0}-3\alpha})}.
 \end{align*}

Analogous to the proof of Theorem \ref{thm:1.1}, we can deduce that
\begin{align*}
& \|u_{L}\|_{X}\leq 3\|u_{0}\|_{\dot{H}^{s_{0}}},\;\|\theta_{L}\|_{Y} \leq 3\|\theta_{0}\|_{\dot{H}^{s_{0}-2\alpha}},\\
&\|L(\theta)\|_{X}\leq C\|P(\theta e_{n})\|_{L^{2}(0,\infty;\dot{H}^{s_{0}-\alpha})}\leq k_{1}\|\theta\|_{Y},\\
&\|\Phi(u,v)\|_{X}\leq  k_{2}\|u\|_{X}\|v\|_{X},
\end{align*}
and
\begin{align*}
\|\Psi(u,\theta)\|_{Y}&\leq C\|u\cdot\nabla \theta\|_{L^{2}(0,\infty;\dot{H}^{s_{0}-3\alpha})}\\&\leq k_{3}\sup_{t\geq0}\|u(t)\|_{\dot{H}^{s_{0}}}\|\theta\|_{L^{2}(0,\infty;\dot{H}^{s_{0}-\alpha})}\\
&\leq
k_{3}\|u\|_{X}\|\theta\|_{Y}.
\end{align*}
Based on these renewed estimates, we can also conclude that if
\begin{align*}
  \|u_{0}\|_{\dot{H}^{s_{0}}}+2k_{1}\|\theta_{0}\|_{\dot{H}^{s_{0}-2\alpha}}\leq\frac{1}{96(k_{2}+k_{3})},
\end{align*}
then
\begin{align*}
  K_{0}:=\|u_{L}\|_{X}+2k_{1}\|\theta_{L}\|_{Y}\leq\frac{1}{32(k_{2}+k_{3})},
\end{align*}
consequently the abstract  equations system (\ref{eqn:bsn2})
has a unique solution $(u,\theta)\in X\times Y$, which is exact the strong solution of (\ref{eqn:bsn1}) on the whole interval $[0,\infty)$, where $\dot{H}^{s_{0}-2\alpha}$ is replaced by $\dot{H}^{s_{0}-3\alpha}$. Besides, similar to  (\ref{est:sol}), we obtain
\begin{align*}
\|u\|_{X}+2k_{1}\|\theta\|_{Y}\leq4K_{0}\leq12(\|u_{0}\|_{\dot{H}^{s_{0}}}+2k_{1}\|\theta_{0}\|_{\dot{H}^{s_{0}-2\alpha}}).
\end{align*}
Moreover, the function $\nabla\pi=Q(\theta e_{n}-u\cdot\nabla u)\in L^{2}(0,\infty;\dot{H}^{s_{0}-\alpha})$ such that $(u,\nabla\pi,\theta)$ fulfills both the equations of (\ref{eqn:bsn3}) in $\dot{H}^{s_{0}-\alpha}$ and $\dot{H}^{s_{0}-3\alpha}$ for a.e. $t\in[0,\infty)$, and
\begin{align*}
   \|\nabla\pi\|_{L^{2}(0,\infty;\dot{H}^{s_{0}-\alpha})}\leq12(\|u_{0}\|_{\dot{H}^{s_{0}}}+2k_{1}\|\theta_{0}\|_{\dot{H}^{s_{0}-\alpha}}).
 \end{align*}
Thus existence part the proof has been completed.

Uniqueness part of the proof is much similar that in the proof of Theorem \ref{thm:1.1}, with the minor differences
\begin{align*}
\|\delta \theta\|_{Y_{t}}\leq k_{3}\|\delta u\|_{X_{t}}\|\theta_{2}\|_{L^{2}(0,t;\dot{H}^{s_{0}-\alpha})}+C\big\|\|u_{1}\|_{\dot{H}^{s_{0}+\varepsilon}}\|\delta\theta\|_{\dot{H}^{s_{0}-\alpha-\varepsilon}}\big\|_{L^{2}(0,t)},
\end{align*}
and
 \begin{align*}
&\quad\big\|\|u_{1}\|_{\dot{H}^{s_{0}+\varepsilon}}\|\delta\theta\|_{\dot{H}^{s_{0}-\alpha-\varepsilon}}\big\|_{L^{2}(0,t)}\\
&\leq\big\|\|u_{1}\|_{\dot{H}^{s_{0}}}^{1-\varepsilon/\alpha}\|u_{1}\|_{\dot{H}^{s_{0}+\alpha}}^{\varepsilon/\alpha}\|\delta \theta\|_{\dot{H}^{s_{0}-2\alpha}}^{\varepsilon/\alpha}\|\delta\theta\|_{\dot{H}^{s_{0}-\alpha}}^{1-\varepsilon/\alpha}\big\|_{L^{2}(0,t)}\\
&\leq\|u_{1}\|_{X}^{1-\varepsilon/\alpha}\|u_{1}\|_{L^{2}(0,t;\dot{H}^{s_{0}+\alpha})}^{\varepsilon/\alpha}\|\delta \theta\|_{Y_{t}}.
\end{align*}
Other inequalities and derivations are omitted here.
\hfill $\Box$

\begin{rem}
With slight modification on the proof of Theorem \ref{thm:1.1}, we can also derive that under the assumption $\frac{1}{2}<\alpha<\frac{2+n}{4}$,  if $u_{0}\in \dot{H}_{\sigma}^{s_{0}}$ with $\|u_{0}\|_{\dot{H}^{s_{0}}}\leq\varepsilon$ for some small number $\varepsilon=\varepsilon(n,\alpha)>0$,
 the initial value problem of the Navier-Stokes equation
  \begin{align*}
\left\{\begin{array}{l}
\partial_{t}u+(-\Delta)^{\alpha}u+u\cdot\nabla u+\nabla \pi=0,\;\;t>0,\;x\in\mathbb{R}^{n};\\
\nabla\cdot u=0,\;\;t>0,\;x\in\mathbb{R}^{n};\\
u(0,x)=u_{0}(x),\;\;x\in\mathbb{R}^{n}.
\end{array} \right.
\end{align*}
  has a unique strong solution $u\in X$ verifying $\|u\|_{X}\leq C\|u_{0}\|_{\dot{H}^{s_{0}}}$
for some the constant $C>0$ independent of $u_{0}$. This is exactly the result derived in \cite{du2018} (where restriction $\alpha>1/2$ should be supplied).  In this sense, results obtained here can be viewed as the extension of \cite{du2018} to the Boussinesq system. Considering that here both the viscosity numbers $\mu,\nu$ are not null, our investigations are also useful supplements to the work \cite{jy2016}.
\end{rem}

\noindent\thanks{{\bf Declarations}}

\noindent\thanks{{\bf Competing interests}
The authors declare that there are no conflicts of interest regarding the publication of
this paper.}

\noindent\thanks{{\bf Authors' contributions} Qinghua Zhang: Conceptualization, Methodology and  Formal analysis. Huiyang Zhang: Writing-Original draft preparation.  Shuokai Yan: Writing-Reviewing and Editing. All authors reviewed and supported the final version of the manuscript.}

\noindent\thanks{{\bf Funding} This research did not receive any specific grant from funding agencies in the public, commercial, or not-for-profit sectors.}


\end{document}